\documentclass[a4paper]{article}
\usepackage{amsmath,amssymb,amsthm}
\newtheorem{theorem}{Theorem}
\newtheorem{lemma}[theorem]{Lemma}
\newtheorem*{corollary}{Corollary}

\theoremstyle{definition}

\newtheorem*{remark}{Remark}
\newtheorem*{example}{Example}

\def\abs#1{\left|#1\right|}
\def\C{\mathbb C}
\def\Cnt#1{{\mathcal C}^{#1}}
\def\csubset{\subset\subset}
\def\eps{\varepsilon}
\def\eqnref#1{(\ref{#1})}
\def\Fin{\mathrm{Fin}}
\def\inw#1{{}^{\circ}\left(#1\right)}
\def\lapprox{\lessapprox}
\def\N{\mathbb N}
\def\ns{\mathrm{ns}}
\def\st{\mathrm{st}}
\def\ster#1{{{}^* #1}}
\def\supp{\mathrm{supp}}
\def\test{\mathcal D}
\def\R{\mathbb R}
\def\Z{\mathbb Z}

\begin{document}
\title{The Local Structure of\\Nonstandard Representatives of Distributions}

\author{Hans Vernaeve\\
Ghent University\\
Dept.\ of Pure Mathematics and Computer Algebra\\
Galg\-laan 2, B-9000 Gent (Belgium)\\
E-mail: {\tt hvernaev@cage.ugent.be}}

\date{}
\maketitle

\begin{abstract}
It is shown that the nonstandard representatives of Schwartz-distri\-butions,
as introduced by K.~D.~Stroyan and W.~A.~J.~Luxemburg in their book
{\it Introduction to the theory of infinitesimals}~\cite{SenL}, are locally
equal to a finite-order derivative of a finite-valued and S-continuous function.
By `equality', we mean a pointwise equality, not an equality in a distributional
sense. This proves a conjecture by M.~Oberguggenberger in~[Z.\ Anal.\ Anwend.\
10 (1991), 263--264]. 
Moreover, the representatives of the zero-distribution are locally equal
to a finite-order derivative of a function assuming only infinitesimal values.
These results also unify the nonstandard theory of distributions by K.~D.~Stroyan and
W.~A.~J.~Luxemburg with the theory by R.~F.~Hoskins and J.~Sousa Pinto
in~[Portugaliae Mathematica 48(2), 195--216].
\end{abstract}



\section{Stroyan and Luxemburg's theory of distributions}
In~\cite[\S 10.4]{SenL}, K.~D.~Stroyan and W.~A.~J.~Luxemburg introduced their nonstandard theory of Schwartz distributions. We give a brief account of the definitions and properties in this theory needed in the sequel. The notations in this section will be used throughout the whole paper (some are different from Stroyan and Luxemburg's). The nonstandard language used is Robinson's.

We will often identify a standard entity $A$ with its image $^\sigma A:=\{\ster x: x\in A\}$ when no confusion is possible.

Let $\Omega$ be an open subset of $\R^n$. Let $\Cnt\infty(\Omega)$ be the 
space of all $\Omega\to\C$-functions possessing continuous derivatives of any order. Let $\test(\Omega)$ be the space of all test-functions on $\Omega$, i.e., all $\Cnt\infty(\Omega)$-functions with compact support contained in $\Omega$ and $\test'(\Omega)$ the space of Schwartz distributions, i.e., continuous linear functionals on $\test(\Omega)$. By $\ns(\ster\Omega)$, we denote the set $\{x\in\ster\Omega: \exists y\in\Omega: x\approx y\}$ of near-standard points of $\ster\Omega$. By $\Fin(\ster\C)$, we denote the set of finite elements of $\ster\C$. By $\st$ we denote the standard part map.

A topological structure is introduced on $\ster\test(\Omega)$ in the following way. We denote by $\partial^\alpha$ the partial derivative of order $\alpha\in\N^n$. A function $\phi\in\ster\test(\Omega)$ is called a \emph{finite} element of $\ster\test(\Omega)$ iff its support is contained in $\ns(\ster\Omega)$ and if $\partial^\alpha \phi(x)\in\Fin(\ster\C)$, for all (finite) multi-indices $\alpha\in\N$ and all $x\in\ster\Omega$. The set of all finite elements of $\ster\test(\Omega)$ will be denoted by $\Fin(\ster\test(\Omega))$.\\
Similarly, $\phi\in\ster\test(\Omega)$ is called an \emph{infinitesimal} element of $\ster\test(\Omega)$ iff its support is contained in $\ns(\ster\Omega)$ and if $\partial^\alpha \phi(x)\approx 0$, for all (finite) multi-indices $\alpha\in\N$ and all $x\in\ster\Omega$. We will write $\phi\approx_{\test} 0$ in this case.

A $\ster{\Cnt\infty}(\Omega)$-function $f$ is called a representative of $T\in\test'(\Omega)$ iff for each $\phi\in\Fin(\ster\test(\Omega))$,
\[\int_{\ster\Omega}f\phi\approx (\ster T)(\phi).\]

It can be shown that every function $f$ in the set
\[D'(\Omega):=\big\{f\in\ster{\Cnt\infty}(\Omega): \int_{\ster\Omega}f\phi\in\Fin(\ster\C),\quad \forall\phi\in\Fin(\ster\test(\Omega))\big\}\]
is a representative of a distribution $T$ by means of the definition
$T(\phi):=\st \int_{\ster\Omega}f\phi$. This unique distribution is called the standard part of $f$ and is denoted by $\st f$.\\
Vice versa, it can be shown that every distribution has a representative in $D'(\Omega)$.

$T\in\ster\test'(\Omega)$ is called S-continuous iff
\begin{equation}\label{eq1}
(\forall\phi\in\ster\test(\Omega))(\phi\approx_\test 0\implies T(\phi)\approx 0).
\end{equation}
It can be shown that every $f\in D'(\Omega)$ is S-continuous as an element of $\ster\test'(\Omega)$. Stroyan and Luxemburg call the elements of $D'(\Omega)$ finite distributions. To avoid the suggestion that $D'(\Omega)$ should be a subset of the space of distributions, and because of the S-continuity as an element of $\ster\test'(\Omega)$, we will call them \emph{S-distributions} instead.

\begin{remark}
A function $f$: $\ster\Omega\to\ster\C$ is called S-continuous iff
\[
x\approx y\implies f(x)\approx f(y),\qquad\forall x,y\in\ster\Omega.
\]
To avoid confusion for elements of $D'(\Omega)$, we will refer to the S-continuity in the sense of eq.~\eqnref{eq1} explicitly as `S-continuity as a linear functional'.
\end{remark}

Two elements $f$, $g$ of $D'(\Omega)$ represent the same distribution iff
\[\int_{\ster\Omega}f\phi\approx \int_{\ster\Omega}g\phi,\quad \forall\phi\in\Fin(\ster\test(\Omega)).\]
In such case, $f$ and $g$ are called $\test'$-infinitely close, and we write
$f\approx_{\test'(\Omega)} g$. If $\Omega$ is fixed in the context and no confusion can exist, we often shortly write $f\approx_{\test'} g$.

\section{The order of an S-distribution}
As it will play a crucial role in proving our results, we recall a theorem about S-continuity which is proved implicitly in~\cite{SenL} (i.e., there is a general theorem on S-continuity from which this theorem follows partly). Also in the context of Banach spaces, characterizations for S-continuity for internal linear maps are well-known (see e.g.~\cite{Wolff}).\\
We write $K\csubset\Omega$ if $K$ is a compact subset of $\Omega$.
\begin{theorem}\label{Scont}
Let $T\in\ster\test'(\Omega)$. Then the following are equivalent:
\begin{enumerate}
\item $T$ is S-continuous
\item $(\forall\phi\in\ster\test(\Omega))$ $(\phi\approx_\test 0 \implies T(\phi)\in\Fin(\ster\C))$
\item $(\forall\phi\in\Fin(\ster\test(\Omega)))$ $(T(\phi)\in\Fin(\ster\C))$
\item $(\forall K\csubset\Omega)$ $(\exists C\in\R)$ $(\exists m\in\N)$
$(\forall \phi\in\ster\test(K))$
\[\big(\abs{T(\phi)}\le C\max_{\abs{\alpha}\le m}
\sup_{x\in\ster K}\abs{\partial^\alpha \phi(x)}\big)\]
\item $(\forall K\csubset\Omega)$ $(\forall \eps\in\R^+)$ $(\exists\delta\in\R^+)$
$(\exists m\in\N)$ $(\forall\phi\in\ster\test(K))$
\[\big(\max\limits_{\abs{\alpha}\le  m}\sup\limits_{x\in
\ster K}\abs{\partial^\alpha\phi(x)}<\delta\implies\abs{T(\phi)}<\eps\big).\]
\end{enumerate}
\end{theorem}
\begin{proof}
$1\Rightarrow 2\Rightarrow 3$: follows using the fact that
$\eps\phi\approx_\test 0$, $\forall \eps\in\ster\R$ with $\eps\approx 0$ and $\forall \phi\in\Fin(\ster\test(\Omega))$.\\
$3\Rightarrow 4$: let $K\csubset\Omega$. Let $m\in\ster\N\setminus\N$ and
$\phi\in\ster\test(K)$. Let 
\[M:=\max\limits_{\abs{\alpha}\le  m}\sup\limits_{x\in
\ster K}\abs{\partial^\alpha\phi(x)}.\]
If $M\ne 0$, $\frac{1}{M}\phi\in\Fin(\ster\test(\Omega))$.
So $\abs{T(\phi)}=M \underbrace{\abs{T(\phi/M)}}_{\in\Fin(\ster\R)}$, and
the internal set
\[\big\{m\in\ster\N: (\forall\phi\in\ster\test(K))
\big(\abs{T(\phi)}\le m \max\limits_{\abs{\alpha}\le  m}\sup\limits_{x\in
\ster K}\abs{\partial^\alpha\phi(x)}\big)\big\}\]
contains all infinite $m$. By underspill, property 4 holds.\\
$4\Rightarrow 5\Rightarrow 1$: follows using the fact that for each $\phi\in\Fin(\ster\test(\Omega))$, there exists $K\csubset\Omega$ such that $\supp\phi\subseteq\ster K$.
\end{proof}

Following Stroyan and Luxemburg, we introduce the notion of S-distri\-bu\-tions 
of finite order.\\
An S-distribution $f$ is of \emph{order} at most $m\in\N$ on $K\csubset\Omega$ iff
\[(\exists C\in\R^+) (\forall \phi\in\ster\test(K)) \Big(\abs{\int_{\ster
\Omega} f \phi}\le C \max_{\abs{\alpha}\le m}\sup_{x\in\ster K}
\abs{\partial^\alpha\phi(x)}\Big)\]
or, equivalently, iff
\[(\exists C\in\R^+) (\forall \phi\in\Fin(\ster\test(K))) \Big(\abs{\int_{\ster
\Omega} f \phi}\le C \max_{\abs{\alpha}\le m}\sup_{x\in\ster K}
\abs{\partial^\alpha\phi(x)}\Big).\]
The smallest $m\in\N$ for which $f$ is of order at most $m$ is (logically)
called the order of $f$.\\
The equivalence of both definitions follows from the fact that for each $\phi\in\ster\test(K)$,
there exists $M\in\ster\R^+$ such that $\phi/M\in\Fin(\ster\test(K))$ (see
the proof of theorem~\ref{Scont}).\\
Any S-distribution $f$ is of some finite order on any given $K\csubset\Omega$.
This follows from theorem~\ref{Scont} applied to the `regular' functional $\phi\mapsto \int_{\ster\Omega} f\phi \in\ster\test'(\Omega)$.

\section{Introduction to the new results in this paper}
In their short section on distributions (which they call a `sketch' themselves), Stroyan and Luxemburg only mention S-distributions of finite order for proving the theorem that every distribution is locally a finite order derivative of a continuous function, by means of the fact (mentioned as an exercise) that any S-distribution of finite order is $\test'$-infinitely close to a finite-order derivative of an S-continuous function $\in D'(\Omega)$. We will show that the order of an S-distribution $f$ is \emph{not} equal to the order of the distribution $\st f$. The difference between these two orders will be the key to give (at least partially) an answer the following questions.\\
What do S-distributions look like? Except from their definition, what are qualitative ways in which they differ from ordinary functions in $\ster{\Cnt\infty}(\Omega)$?\\
How much can two representatives of the same distribution differ? Except from
the fact that they are $\test'$-infinitely close, are there qualitative ways in
which this difference can be described?\\
It may be clear from the following example that there is hardly any pointwise way in which different representations from a given distribution coincide in general.

\begin{example} \label{nopointrepr}
For each $k\in\Z$ and $\omega\in\ster\N\setminus\N$, the function
$\omega^k\sin(\omega x)\in\ster{\Cnt\infty}(\R)$ is a representative of the zero-distribution ($\in\test'(\R)$).
\end{example}
\begin{proof}
For $k<0$, $f_k(x)=\omega^k\sin(\omega x)\approx 0$, $\forall x\in\ster\R$, so
$f_k\approx_{\test'} 0$. As it is well-known that the distributional derivatives coincide with the derivatives of the representatives, also the second derivative $f_k''=-f_{k+2}\approx_{\test'} 0$. Inductively,
$f_k\approx_{\test'} 0$, $\forall k\in\N$.
\end{proof}
In the example, the method to find heavily irregular representatives of the zero-distribution was by taking derivatives of a function that assumes
infinitesimal values. We will prove that no other irregularities can exist, i.e., that every $f\approx_{\test'} 0$ is
(locally) pointwise equal to some finite order derivative of a
$\ster{\Cnt\infty}(\Omega)$-function assuming only infinitesimal values.\\
Similarly, we will prove that every $f\in D'(\Omega)$ is
(locally) pointwise equal to some finite order derivative of an S-continuous and finite-valued $\ster{\Cnt\infty}(\Omega)$-function.\\
The last of these two assertions was already mentioned (for $\Omega=\R^n$ and
omitting the S-continuity) in \cite[Prop.~2.10]{Ober-ns} in the nonstandard
language of Nelson, but, as it appears from the correction to \cite{Ober-ns}, it
still remained unproved.\\
Although such theorems are of a fashion similar to the classical local
representation theorem of distributions, the distributional order cannot be a
measure for the order of the derivative in our representation theorems: already
for the zero-distribution, which is trivially of order $0$, the order of the
derivative may be arbitrary large.

\section{Proofs of the new results}
First, we point out more explicitly that the order of an S-distribution is not equal to the distributional order of its standard part. For $x$, $y$ $\in\ster\R$, we write $x\lapprox y$ iff $x<y$ or $x\approx y$.
\begin{theorem}
Let $f\in D'(\Omega)$ and $K\csubset\Omega$. Then the (distributional) order of $\st f$ on $K$ is the smallest $m\in\N$ such that
\begin{equation}\label{distr-orde}
(\exists C\in\R^+) (\forall \phi\in\Fin(\ster\test(K))) \Big(\abs{\int_{\ster
\Omega} f \phi}\lapprox C \max_{\abs{\alpha}\le m}\sup_{x\in\ster K}
\abs{\partial^\alpha\phi(x)}\Big).
\end{equation}
\end{theorem}
\begin{proof}
1. Let the order of $T:=\st f$ on $K$ be at most $m$, i.e.\ (by transfer),
\[(\exists C\in\R^+) (\forall \phi\in\ster\test(K)) (\abs{\ster T(\phi)}\le C \max_{\abs{\alpha}\le m}\sup_{x\in\ster K}
\abs{\partial^\alpha\phi(x)}).\]
Since for $\phi\in\Fin(\ster\test(\Omega))$, $\ster T(\phi)\approx
\int_{\ster\Omega} f\phi$, we find that formula~\eqnref{distr-orde} holds for
this $m$.\\
2. On the other hand, suppose that formula~\eqnref{distr-orde} holds for some
$m\in\N$. Again by the fact that for $\phi\in\Fin(\ster\test(\Omega))$, $\ster T(\phi)\approx
\int_{\ster\Omega} f\phi$ (with $T=\st f$), we have in particular that
\[(\exists C\in\R^+) (\forall \phi\in\test(K)) (\abs{\ster
T(\ster\phi)}\lapprox C \max_{\abs{\alpha}\le m}\sup_{x\in\ster K}
\abs{\partial^\alpha\ster\phi(x)}).\]
Since both sides of the $\lapprox$-inequality are standard numbers, we actually
have a $\le$-inequality, and the (distributional) order of $T$ on $K$ is at most $m$.
\end{proof}
\begin{corollary}
The order of an S-distribution $f$ is at least the distributional order of $\st f$.
\end{corollary}

The following example shows that the difference of the two orders can be arbitrary large.
\begin{example}\label{order-ex}
Consider $f(x)=\omega^k \sin(\omega x)$, with $\omega\in\ster\N\setminus\N$.
It has order $k$ on every compact $K\csubset\R$. On the other hand, $f\approx_{\test'} 0$ (see example \ref{nopointrepr}),
so the order of the corresponding standard distribution is $0$.
\begin{proof}
Let $\phi\in\ster\test(K)$. For some $R\in\R$, $K\subseteq [-R,R]$. Then by
partial integration,
\[\int_{\ster\R}f\phi=(-1)^k\int_{\ster\R}g(x)\phi^{(k)}(x)\,dx\]
with $g^{(k)}=f$, so we can choose $g(x)\in\{\pm\sin(\omega x),\pm\cos(\omega x)\}$.
So
\[\abs{\int_{\ster\R}f\phi}
\le 2R \sup_{x\in\ster K}\abs{g(x)}\sup_{x\in\ster K}\abs{\phi^{(k)}(x)}
\le 2R\sup_{x\in\ster K}\abs{\phi^{(k)}(x)},\]
so the order is at most $k$.\\
To see that the order is at least $k$, let
$\phi_0\in\test(K)$ with $\int\phi_0=1$ and let
$\phi(x):= \sin(\omega x)\phi_0(x)$. Then
\[\frac{1}{\omega^k}\int_{\ster\R}f\phi=\frac{1}{2}\int_{\ster\R}(1-\cos(2\omega
x))\phi_0(x)\,dx\approx \frac{1}{2},\]
since $\cos(2\omega x)\approx_{\test'} 0$
(similarly as in example~\ref{nopointrepr}). On the other hand, for each
$j\in\N$, $\sup_{x\in\ster K}\abs{\phi^{(j)}(x)}\le M \omega^{j}$ for some
$M\in\R$, so for this $\phi\in\ster\test(K)$, $\abs{\int_{\ster\R}f\phi}>
C\max_{j\le k-1}\sup_{x\in\ster K}\abs{\phi^{(j)}(x)}$, $\forall C\in\R$.
\end{proof}
\end{example}

Next, we will prepare our main results. First, we show that distributional
anti-derivatives can be dealt with on representatives. To our know\-ledge, such
a theorem is not available in the nonstandard literature. Just for convenience,
we only deal with partial derivatives in the first variable.\\
We introduce the
following notation: for $x=(x_1,\ldots,x_n)\in\ster\R^n$, we will write $\tilde
x_i:= (x_1,\ldots,x_{i-1},x_{i+1},\ldots,x_n)$. Similarly, for $i<j$, $\tilde
x_{i,j}$ $=$ $(x_1$,\dots, $x_{i-1}$, $x_{i+1}$,\dots, $x_{j-1}$,
$x_{j+1}$,\dots, $x_n)$ and so on for $\tilde x_{i,j,k}$, \dots

\begin{lemma}\label{anti-derivative}
Let $\Omega$ be an open interval (i.e., it is the Cartesian product of $n$
one-dim.\ intervals). Let $T\in \test'(\Omega)$ and $f$ a representative of $T$. Then there exists an S-distribution $g\in D'(\Omega)$ with $\partial_1 g=f$.
As a consequence, $g$ determines a distribution $U$ with $\partial_1 U= T$.
\end{lemma}
\begin{proof}
1. In order to get some insight in the proof, we first consider the one-dimensional case.\\
Choose $F\in\ster{\Cnt\infty}(\Omega)$ such that $F'=f$ on $\Omega$. We can only expect $F$ to be an S-distribution if the
integration constant is well-chosen. So, we seek $C\in\ster\C$ such that
$\int_{\ster\R} (F+C)\phi\in\Fin(\ster\C)$, $\forall\phi\in\Fin(\ster\test(\Omega))$. Now fix
$\phi_0\in\test(\Omega)$, with $\int_\R\phi_0=1$. Then the previous condition specifies
to $\int_{\ster\R} F\ster\phi_0+C\in\Fin(\ster\C)$. As a finite change in the constant
doesn't influence the S-distributional character of $F+C$, we can put
$C:=-\int_{\ster\R} F\ster\phi_0$. Then, for any $\phi\in\Fin(\ster\test(\Omega))$,
\[\int_{\ster\R} (F+C)\phi=\int_{\ster\R} F(t)
\big(\underbrace{\phi(t)-\Big(\int_{\ster\R}\phi\Big)\ster\phi_0(t)}_{=:\psi(t)\in\Fin(\ster\test(\Omega))}\big)
\,dt.\]
As $\int_\ster\R\psi=0$, $\psi^{(-1)}(x):=\int_{-\infty}^x \psi\in\Fin(\ster\test(\Omega))$, and
by partial integration,
\[\int_{\ster\R} (F+C)\phi=-\int_{\ster\R} f\psi^{(-1)}\in\Fin(\ster\C),\]
since $f$ is an S-distribution.\\
2. In the general case, we choose an arbitrary anti-derivative $F$ of $f$ in the
first variable (on $\Omega$). E.g., if
$\Omega=(a_1,b_1)\times\cdots\times(a_n,b_n)$ ($a_i$,
$b_i\in\R\cup\{-\infty,+\infty\}$), then for any $a_1<c<b_1$,
$\int_{c}^{x_1}f(t,\tilde x_1)\,dt$ is such an anti-derivative). An
anti-derivative is determined up to a function $G(\tilde x_1)$. Now it turns
out that, for a fixed $\phi_0\in\test((a_1,b_1))$ with $\int_\R \phi_0=1$,
$G(\tilde x_1)=-\int_{\ster\R} F(t,\tilde x_1)\ster\phi_0(t)\,dt$ is a good choice:
for any $\phi\in\Fin(\ster\test(\Omega))$,
\begin{multline*}
\int_{\ster\R^n}(F(x)+G(\tilde x_1))\phi(x)\,dx\\
=\int_{\ster\R^n}F(x)
\big(\underbrace{\phi(x)-\Big(\int_{\ster\R}\phi(u,\tilde
x_1)\,du\Big)\ster\phi_0(x_1)}_{=:\psi(x)}
\big)\,dx.
\end{multline*}
As $\Omega$ is an interval, $\psi\in\Fin(\ster\test(\Omega))$. Moreover,
$\int_{\ster\R}\psi(t,\tilde x_1)\,dt=0$, $\forall\tilde x_1\in\ster\R^{n-1}$, so
$\chi(x):=\int_{-\infty}^{x_1}\psi(t,\tilde x)\,dt\in\Fin(\ster\test(\Omega))$ and similarly as in the one-dimensional case, we find that
$\int_{\ster\R^n}(F(x)+G(\tilde x_1))\phi(x)\,dx\in\Fin(\ster\C)$.
\end{proof}

\begin{lemma}\label{loc-rep2}
Let $f\in D'(\Omega)$ of order $\le m$ on an interval $K\csubset\Omega$, $m>0$.
Then there exists $g\in D'(\Omega)$ of order $\le m-1$ on $K$
such that $\partial_1\cdots\partial_n g=f$ on $\ster K$.
\end{lemma}
\begin{proof}
Let $K=[a_1,b_1]\times\cdots\times[a_n,b_n]$.
We will show that, if $f$ satisfies
\[\abs{\int_{\ster\Omega} f\phi}\le C \sup_{x\in\ster
K}\abs{\partial^{(k,\alpha)}
\phi(x)}, \quad\forall\phi\in\Fin(\ster\test(K))\]
for some $C\in\R$, $k\in\N$ and $\alpha\in\N^{n-1}$,
then the anti-derivative $g(x)=F(x)+G(\tilde x_1)$ in the first variable
defined in lemma~\ref{anti-derivative} satisfies
\[\abs{\int_{\ster\Omega} g\phi}\le C' \max_{j\le l}\sup_{x\in\ster K}
\abs{\partial^{(j,\alpha)} \phi(x)}, \quad\forall\phi\in\Fin(\ster\test(K))\]
with $C'\in\R$ and $l=\max(k-1,0)$.\\
Let $\phi\in\Fin(\ster\test(K))$. With
$\psi$, $\chi\in\Fin(\ster\test(K))$ as in lemma~\ref{anti-derivative},
we have
\[
\abs{\int_{\ster\Omega} g\phi}
=\abs{\int_{\ster\Omega} f \chi}
\le C \sup_{x\in\ster K}\abs{\partial^{(k,\alpha)}
\chi(x)}
= C \sup_{x\in\ster K}\abs{\partial^{(0,\alpha)}\partial_1^{k}
\partial_1^{-1}\psi(x)}.
\]
In case $k=0$, we have for $x\in\ster K$ that
\[\abs{\partial^{(0,\alpha)}\partial_1^{-1}\psi(x)}
=\abs{\int_{-\infty}^{x_1}\partial^{(0,\alpha)}\psi(t_1,\tilde x_1)\,dt_1}
\le(b_1-a_1)\sup_{x\in\ster K}\abs{\partial^{(0,\alpha)}\psi(x)},\]
so in any case we have (for some $C'$, $C''$ $\in\R$, independent of $\phi$)
\begin{align*}
\abs{\int_{\ster\Omega} g\phi}
&\le C'\sup_{x\in\ster K}\abs{\partial^{(l,\alpha)}\phi(x)}
+C'\sup_{x\in\ster K}
\abs{D^l \,\ster\phi_0(x_1)\int_{\ster\R}\partial^{(0,\alpha)}\phi(u,\tilde x_1)\,du}\\
&\le C''\max_{j\le l}\sup_{x\in\ster K}\abs{\partial^{(j,\alpha)}\phi(x)}.
\end{align*}
Since $g$ is well-defined on $\ster\Omega'$, for some interval $\Omega'\subseteq\Omega$ with $K\csubset\Omega'$, we can use $\phi_0\in\test(\Omega')$ with $\phi_0=1$ on $K$ to ensure that $g\ster\phi_0\in D'(\Omega)$ without changing the values on $\ster K$.\\
If we repeatedly apply also the analogous result for the variables $x_2$, \ldots,
$x_n$, we finally conclude that the order of the primitive
$(\partial_1\cdots\partial_n)^{-1} f$ has decreased (if $m>0$).
\end{proof}

For $K\csubset\Omega$, we call $L^{\infty}(K)$ the space of all (standard)
bounded and (Lebesgue-)measurable functions $f$: $\Omega\to\C$ with support contained in $K$.
\begin{lemma}\label{order-zero}
Let $K\csubset\Omega$ an interval. An S-distribution $f$ is of order zero on
$K$ iff
\[(\exists C\in\R^+) (\forall \phi\in\ster L^{\infty}(K)) \Big(\abs{\int_{\ster
\Omega} f \phi}\le C \sup_{x\in\ster K} \abs{\phi(x)}\Big).\]
\end{lemma}
\begin{proof}
Let $f\in\Cnt{\infty}(\Omega)$ and $\phi\in
L^{\infty}(K)$. Then by a classical density theorem,
it is clear that there exists some $h\in\test(K)$ such that
\[\abs{\int_{\Omega} f\phi-\int_{\Omega} f h}\le\sup_{x\in K} \abs{\phi(x)}\quad \&\quad \sup_{x\in K}\abs{h(x)} \le 2\sup_{x\in K}\abs{\phi(x)}.\]
By transfer, we have
$(\forall f\in\ster{\Cnt\infty}(\Omega))$ $(\forall\phi\in \ster L^\infty(K))$
$(\exists h\in\ster\test(K))$
\[\Big(\abs{\int_{\ster \Omega} f\phi-\int_{\ster\Omega}
fh}\le\sup_{x\in\ster K}\abs{\phi(x)} \quad \&\quad
\sup_{x\in\ster K}\abs{h(x)}\le 2\sup_{x\in\ster K}\abs{\phi(x)}\Big).\]
If in particular $f$ is an S-distribution of order $0$ on $K$, then
\[(\exists C\in\R^+) (\forall h\in\ster\test(K))\big(\abs{\int_{\ster\Omega}
fh}\le C\sup_{x\in\ster K}\abs{h(x)}\big).\]
The result follows by combining these two formulas.
\end{proof}

\begin{lemma}\label{loc-rep3}
Let $f\in D'(\Omega)$. Suppose that $f$ is of order zero on a (standard)
interval $K=[a_1,b_1]\times\cdots\times [a_n,b_n]\csubset\Omega$. Then
\begin{enumerate}
\item there exists $g\in\ster{\Cnt\infty}(\Omega)$ which is bounded on $\ster K$
by a standard constant and such that $\partial_1\cdots\partial_n g=f$ on $\ster K$.
\item there exists $h\in\ster{\Cnt\infty}(\Omega)$ which is S-continuous and
bounded by a standard constant on $\ster K$ and such that
$\partial^2_1\cdots\partial^2_n h=f$ on $\ster K$.
\end{enumerate}
\end{lemma}
\begin{proof}
1. Let $x=(x_1,\ldots,x_n)$ and $t=(t_1,\ldots,t_n)$. For $A\subset \Omega$, we
denote the characteristic function of $A$ by $\chi_A$. Then (for $x\in\ster K$)
\[g(x):=\int_{a_1}^{x_1}dt_1\cdots\int_{a_n}^{x_n} f(t)\, dt_n=\int_{\ster\Omega}
f\chi_{[a_1,x_1]\times\cdots\times[a_n,x_n]}\] clearly satisfies
$\partial_1\cdots\partial_n g=f$ on $\ster K$. Further, applying the previous lemma with
$\phi=\chi_{[a_1,x_1]\times\cdots\times[a_n,x_n]}\in\ster L^\infty(K)$ (if
$x\in \ster K$), we find $C\in\R^+$ such that
\[(\forall x\in \ster K)\big(\abs{g(x)}\le C
\underbrace{\sup_{x\in\ster K} \abs{\phi(x)}}_{=1}\big).\]
2. If $g$ satisfies the conditions from part~1, then (for $x\in\ster K$)
\[h(x):=\int_{a_1}^{x_1}dt_1\cdots\int_{a_n}^{x_n} g(t)\, dt_n\]
clearly satisfies $\partial^2_1\cdots\partial^2_n h=f$ on $\ster K$. Further,
for $\eps\approx 0$, $\eps>0$,
\[\abs{h(x_1+\eps,\tilde x_1)-h(x)}=
\abs{\int_{x_1}^{x_1+\eps}dt_1\cdots\int_{a_n}^{x_n} g(t)\, dt_n}\le
C\eps\prod_{i\ne 1}(b_i-a_i)\approx 0\]
and similarly for the other variables. So $h(x)\approx h(y)$ as soon as $x\approx
y$ ($x$, $y$ $\in \ster K$). Further, $\abs{h(x)}\le C \prod_{i}(b_i-a_i)\in\Fin(\ster\C)$, $\forall x\in\ster K$.
\end{proof}

We are now ready to prove the first main result.
\begin{theorem}\label{local-distr}
Let $f\in \ster{\Cnt\infty}(\Omega)$. Then $f\in D'(\Omega)$ iff for each $K\csubset\Omega$, there exists a
$g\in D'(\Omega)$ which is finite-valued and S-continuous on $\ster K$
and such that $f$ is a finite order derivative of $g$ on $\ster K$.
\end{theorem}
\begin{proof}
$\Leftarrow$: follows using the fact that for each $\phi\in\Fin(\ster\test(\Omega))$, there exists
$K\csubset\Omega$ such that $\supp\phi\subseteq\ster K$.\\
$\Rightarrow$:
1. We first consider the special case where $K\csubset\Omega$ is an interval.\\
Take an interval $K'\csubset \Omega$ with $K\csubset \inw{K'}$, the (topological) interior of $K'$. Since $f$ has a finite order $m$ on $K'$, we find, by repeatedly applying
lemma~\ref{loc-rep2}, some $\tilde g\in D'(\Omega)$ of order zero on $K'$ such that
$(\partial_1\cdots\partial_n)^m \tilde g=f$ on $\ster K'$. By lemma~\ref{loc-rep3}, we
find $h\in\ster{\Cnt\infty}(\Omega)$ which is finite and S-continuous on $\ster K'$
and such that $(\partial_1\cdots\partial_n)^{m+2} h=f$ on $\ster K'$. If 
$\phi_0\in\test(K')$ with $\phi_0=1$ on $K$, then $g:=h\ster\phi_0\in D'(\Omega)$ has
the required properties.\\
2. We consider the special case where $f(x)=0$, $\forall x\notin
\ns(\ster\Omega)$.\\
Then $f$ can be extended to a $\ster{\Cnt\infty}(\R^n)$-function,
setting $f(x):=0$ if $x\in\ster\R^n\setminus\ns(\ster\Omega)$. We claim that this extension $\in D'(\R^n)$. There exists $K_0\csubset\Omega$ such that $f(x)=0$ outside $\ster K_0$. Choose $\phi_0\in\test(\Omega)$ with $\phi_0=1$ on $K_0$. Then for any $\phi\in\Fin(\ster\test(\R^n))$,
\[\int_{\ster\R^n} f\phi =\int_{\ster\Omega}f\underbrace{\ster\phi_0\phi}_{\in\Fin(\ster\test(\Omega))} \in\Fin(\ster\C).\]
Now let $K\csubset\Omega$ arbitrarily. Since $K\subseteq L\csubset\R^n$, with $L$ an interval (possibly
$L\not\subseteq\Omega$), we conclude from part~1 that there exists a
$g\in D'(\R^n)$ which is finite and S-continuous on $\ster L$
and such that (the extended) $f$ is a finite order derivative of $g$ on $\ster L$.
The restriction of $g$ to $\ster\Omega$ has the required properties.\\
3. In the general case, let $K\csubset\Omega$. Taking $\phi_0\in\test(\Omega)$
with $\phi_0=1$ on $K$, we apply part~2 on $f\ster\phi_0\in D'(\Omega)$.
\end{proof}

The second main result will follow from the previous theorem together with some additional lemmas.
\begin{lemma}
Let $\Omega=(a_1,b_1)\times\cdots(a_n,b_n)\subseteq\R^n$ be an open interval
(possibly $a_i=-\infty$, $b_i=+\infty$). Let
$\tilde\Omega:=(a_2,b_2)\times\cdots(a_n,b_n)\subseteq\R^{n-1}$.
Let $f\in\ster{\Cnt\infty}(\Omega)$ be independent of $x_1$, so it can be
identified with a $\ster{\Cnt\infty}(\tilde\Omega)$-function. Then
\begin{enumerate}
\item $f(\tilde x_1)\in D'(\Omega)\iff f(\tilde x_1)\in D'(\tilde\Omega)$.
\item $f(\tilde x_1)\approx_{\test'(\Omega)} 0\iff f(\tilde x_1)\approx_{\test'(\tilde\Omega)} 0$.
\end{enumerate}
As a consequence, the expression $f(\tilde x_1)\approx_{\test'} 0$ is unambiguous.
\end{lemma}
\begin{proof}
1. $\Rightarrow$: Let $f(\tilde x_1)\in D'(\Omega)$.
Fix $\psi(x_1)\in\Fin(\ster\test(a_1,b_1))$ with $\int_{\ster\R}\psi=1$. Choose
$\phi(\tilde x_1)\in\Fin(\ster\test(\tilde\Omega))$ arbitrarily. Then
$\psi(x_1)\phi(\tilde x_1)\in\Fin(\ster\test(\Omega))$, so
\[
\Fin(\ster\C)\ni\int_{\ster\Omega}f(\tilde x_1)\psi(x_1)\phi(\tilde x_1)\,dx
=\underbrace{\int_{a_1}^{b_1}\psi(x_1)\,dx_1}_{=1}
\int_{\ster{\tilde\Omega}}f(\tilde x_1)\phi(\tilde x_1)\,d\tilde x_1,
\]
which means that $f(\tilde x_1)\in D'(\tilde\Omega)$.\\
$\Leftarrow$: Let $f(\tilde x_1)\in D'(\tilde\Omega)$.
For any $\phi\in\Fin(\ster\test(\Omega))$ and $c\in\ns(\ster(a_1,b_1))$, the map
$\tilde x_1\mapsto \phi(c,\tilde x_1)\in\Fin(\ster\test(\tilde\Omega))$, so
\[\psi(c):=\int_{\ster{\tilde\Omega}}f(\tilde x_1)\phi(c,\tilde x_1)\,d\tilde
x_1\in\Fin(\ster\C).\]
Further, for some $K\csubset(a_1,b_1)$,
if $c$ lies outside $\ster K$, $\psi(c)=0$. So
\[\int_{\ster\Omega}f(\tilde x_1)\phi(x)\,dx=\int_{\ster
K}\psi(x_1)\,dx_1\in\Fin(\ster\C),\]
which means that $f(\tilde x_1)\in D'(\Omega)$.\\
2. Similar.
\end{proof}

The following lemmas could be considered as exercises in distribution theory.
To our knowledge, they are not widely known. Therefore, we will include a nonstandard version with proof.
\begin{lemma}
Let $\Omega$ be an open interval.\\
If $f\in D'(\Omega)$ and $\partial^\alpha f\approx_{\test'} 0$,
then there exist $g_{ij} \in D'(\Omega)$ such that
\[f(x)\approx_{\test'}
\sum_{i=1}^n\sum_{j=0}^{\alpha_i-1} g_{ij}(\tilde x_i)x_i^j.\]
\end{lemma}
\begin{proof}
1. We first show that, if $F\in D'(\Omega)$ and $\partial_1 F\approx_{\test'}0$,
then $F$ is $\test'$-infinitely close to a $D'(\Omega)$-function which doesn't
depend on $x_1$.\\
If we choose $G(\tilde x_1)$ as in lemma~\ref{anti-derivative},
we see that for all $\phi\in\Fin(\ster\test(\Omega))$,
\[\int_{\ster\R^n}(F(x)+G(\tilde x_1))\phi(x)\,dx=\int_{\ster\R^n}(\partial_1
F)(x)\chi(x)\,dx\approx 0\]
with $\chi\in\Fin(\ster\test(\Omega))$ as in lemma~\ref{anti-derivative}.
So $F(x)\approx_{\test'}-G(\tilde x_1)$.\\
2. Now suppose that $f\in D'(\Omega)$ and
\begin{equation}\label{eq3}
\partial_1 f(x)\approx_{\test'}
\sum_{i=1}^n\sum_{j=0}^{m_i}g_{ij}(\tilde x_i)x_i^j
\end{equation}
for some $g_{ij}\in D'(\Omega)$.
We will show that
\[f(x)\approx_{\test'}
\sum_{i=1}^n\sum_{j=0}^{\tilde m_i}\tilde g_{ij}(\tilde x_i)x_i^j\]
for some $\tilde g_{ij}\in D'(\Omega)$, $\tilde m_1=m_1+1$, $\tilde
m_2=m_2$, \dots, $\tilde m_n=m_n$.\\
We notice that the right-hand side of eq.~\eqnref{eq3} is equal to
\[\partial_1 \left(\sum_{j=0}^{m_1} g_{1j}(\tilde x_1)\frac{x_1^{j+1}}{j+1}
+\sum_{i=2}^n\sum_{j=0}^{m_i} (\partial_1^{-1}g_{ij})(\tilde x_i) x_i^j\right).
\]
From the explicit construction of the primitives $\partial_1^{-1}g_{ij}$ in
lemma~\ref{anti-derivative}, it is immediate that also they are independent of
$x_i$. Then applying part~1 on the difference of both sides in eq.~\eqnref{eq3},
we find that there exists $G(\tilde x_1)\in D'(\Omega)$ such that
\[f(x)\approx_{\test'}G(\tilde x_1)+\sum_{j=0}^{m_1} g_{1j}(\tilde x_1)\frac{x_1^{j+1}}{j+1}
+\sum_{i=2}^n\sum_{j=0}^{m_i} (\partial_1^{-1}g_{ij})(\tilde x_i) x_i^j\]
which has the required form.\\
3. Now the theorem follows inductively using part~2 and the analogous formulas for the
other variables ($\ne 1$), also using the fact that if $f\in D'(\Omega)$,
then $\partial^\beta f\in D'(\Omega)$, $\forall\beta\in\N^n$.
\end{proof}

\begin{lemma}\label{Scont-prim}
Let $\Omega$ be an open interval.
Let $f\in \ster{\Cnt\infty}(\Omega)$
be S-continuous and finite-valued on $\ns(\ster\Omega)$ and suppose that
$\partial^\alpha f\approx_{\test'} 0$. Then there exist
$g_{ij}\in \ster{\Cnt\infty}(\Omega)$ which are S-continuous and finite-valued on
$\ns(\ster\Omega)$ such that
\[f(x)\approx_{\test'}
\sum_{i=1}^n\sum_{j=0}^{\alpha_i-1} g_{ij}(\tilde x_i)x_i^j.\]
\end{lemma}
\begin{proof}
First notice that a $\ster{\Cnt\infty}(\Omega)$-function which is finite-valued on $\ns(\ster\Omega)$, is $\in D'(\Omega)$.
Let $\Omega=(a_1,b_1)\times\cdots\times(a_n,b_n)$
and $\tilde\Omega:=(a_2,b_2)\times\cdots\times(a_n,b_n)$.
Let $\partial^\alpha f\approx_{\test'} 0$ and let $\alpha=:(\alpha_1,\tilde\alpha)$, $\tilde\alpha\in\N^{n-1}$.
By the previous lemma,
\begin{equation}\label{apply-lemma}
f(x)\approx_{\test'}\sum_{i=1}^n\sum_{j=0}^{\alpha_i-1} h_{ij}(\tilde x_i)x_i^j,
\end{equation}
with $h_{ij}\in D'(\Omega)$.
Now consider $c\in\ns\ster(a_1,b_1)$ arbitrarily. Fix $\psi(x_1)\in\test(\R)$
with $\int_{\R}\psi=1$ and $\psi\ge 0$. Let $\psi_m(x_1):=m\psi(mx_1)$, $\forall
m\in\ster\N$. Let $\phi(\tilde x_1) \in\Fin(\ster\test(\tilde\Omega))$ arbitrarily. Since
\[\partial^{(0,\tilde\alpha)}f(x)\approx_{\test'}\sum_{j=0}^{\alpha_1-1}
\partial^{\tilde\alpha}h_{1j}(\tilde x_1)x_1^j,\]
we have for sufficiently large $m\in\N$ (such that $\supp(\psi_m(c-x_1))\subset\ns\ster(a_1,b_1)$) that
\begin{equation}\label{eq4}
\int_{\ster\Omega}\partial^{(0,\tilde\alpha)}f(x)\psi_m(c-x_1)\phi(\tilde x_1)\,dx
\approx\sum_{j=0}^{\alpha_1-1} \int_{\ster\Omega}
\partial^{\tilde\alpha}h_{1j}(\tilde x_1)x_1^j\psi_m(c-x_1)\phi(\tilde x_1)\,dx.
\end{equation}
By Robinson's sequential lemma, this also holds for some
$\omega\in\ster\N\setminus\N$. If $\tilde x_1\in\ns(\ster{\tilde\Omega})$, the map $x_1\to f(x)$ is S-continuous on $\ns(\ster(a_1,b_1))$.
Then $\forall \tilde x_1\in\ns(\ster{\tilde\Omega})$,
\begin{align*}
&\abs{\int_{a_1}^{b_1} f(x)\psi_\omega(c-x_1)\,dx_1 - f(c,\tilde
x_1)}\\
=&\abs{\int_{a_1}^{b_1} (f(x)-f(c,\tilde x_1))\psi_\omega(c-x_1)\,dx_1}
\le\sup_{x_1\in \supp\psi_\omega}\abs{f(x)-f(c,\tilde x_1)}\approx 0
\end{align*}
since $\supp\psi_\omega$ contains only infinitesimals and $\int_{\ster\R}\abs{\psi_\omega}=1$.
In particular, they are $\test'$-infinitely close. So also
\[\int_{a_1}^{b_1} \partial^{(0,\tilde\alpha)}f(x)\psi_\omega(c-x_1)\,dx_1
\approx_{\test'} \partial^{\tilde\alpha}f(c,\tilde x_1).\]
On the other hand,
\begin{multline*}
\int_{\ster\Omega} \partial^{\tilde\alpha}h_{1j}(\tilde x_1)x_1^j
\psi_\omega(c-x_1)\phi(\tilde x_1)\,dx\\
=\underbrace{\int_{\ster{\tilde\Omega}} \partial^{\tilde\alpha}h_{1j}(\tilde x_1)
\phi(\tilde x_1)\,d\tilde x_1}_{\in\Fin(\ster\C)}
\underbrace{\int_{a_1}^{b_1} x_1^j\psi_\omega(c-x_1)\,dx_1}_{\approx c^j},
\end{multline*}
so we find from equation \eqnref{eq4} that for each $c\in\ns(\ster(a_1,b_1))$
\[\partial^{\tilde\alpha}f(c,\tilde x_1)\approx_{\test'}
\sum_{j=0}^{\alpha_1-1} c^j \partial^{\tilde\alpha}h_{1j}(\tilde x_1).\]
Now choose $\alpha_1$ different values $c_i\in\ns(\ster(a_1,b_1))$, with
$c_i\not\approx c_j$ if $i\ne j$. Then we find a linear system with $\alpha_1$
equations and $\alpha_1$ unknown functions $\partial^{\tilde\alpha}h_{1j}$.
The determinant of the system is a Vandermonde-determinant equal to
$\prod_{i<j}(c_j-c_i)\not\approx 0$.
Therefore, each $\partial^{\tilde\alpha}h_{1j}(\tilde x_1)$ is
$\test'$-infinitely close to a $\Fin(\ster\C)$-linear combination of the
$\partial^{\tilde\alpha} f(c_j,\tilde x_1)$, which we call $\partial^{\tilde\alpha} g_{1j}(\tilde x_1)$. So $g_{1j}(\tilde x_1)\in\ster{\Cnt\infty}(\Omega)$ are S-continuous and finite-valued on $\ns(\ster\Omega)$. By the previous lemma (applied to $\tilde\Omega\subseteq\R^{n-1}$),
\[h_{1j}(\tilde x_1)\approx_{\test'} g_{1j} (\tilde x_1)+ \sum_{i=2}^n\sum_{k=0}^{\alpha_i-1} \tilde h_{ik}(\tilde x_{1i})x_i^k,\]
for some $\tilde h_{ik}\in D'(\Omega)$.
Substituting these expressions, together with the analogous expressions for $h_{ij}(\tilde x_i)$ (with $i>1$), in formula~\eqnref{apply-lemma} yields that
\begin{equation}\label{induction-hyp}
f(x)\approx_{\test'}\sum_{i=1}^n\sum_{j=0}^{\alpha_i-1} g_{ij}(\tilde x_i)x_i^j
+\sum_{1\le i_1<i_2\le n}\sum_{j_1=0}^{\alpha_{i_1}-1}\sum_{j_2=0}^{\alpha_{i_2}-1} h_{i_1i_2j_1j_2}(\tilde x_{i_1,i_2})x_{i_1}^{j_1}x_{i_2}^{j_2}\,,
\end{equation}
for some $h_{i_1i_2j_1j_2}\in D'(\Omega)$, since multiplication by $x_i$ preserves the $\approx_{\test'}$-equality.

We now proceed inductively and show that
\begin{equation}\label{induction-result}
f(x)\approx_{\test'}\sum_{i=1}^n\sum_{j=0}^{\alpha_i-1} g_{ij}(\tilde x_i)x_i^j
+\!\!\!\sum_{\substack{1\le i_1<i_2\ \ \\ \ \ <i_3\le n}}\!\sum_{j_1,j_2,j_3} h_{i_1i_2i_3j_1j_2j_3}(\tilde x_{i_1,i_2,i_3})x_{i_1}^{j_1}x_{i_2}^{j_2}x_{i_3}^{j_3}\,,
\end{equation}
for some $g_{ij}\in\ster{\Cnt\infty}(\Omega)$, S-continuous and finite-valued on
$\ns(\ster\Omega)$ and some $h_{i_1i_2i_3j_1j_2j_3}\in D'(\Omega)$.\\
The proof is similar. Let $F:=f- \sum_{i=1}^n\sum_{j=0}^{\alpha_i-1} g_{ij}(\tilde x_i)x_i^j$. Let $\alpha=:(\alpha_1,\alpha_2,\tilde\alpha)$, $\tilde\alpha\in\N^{n-2}$. Let $\tilde\Omega:=(a_3,b_3)\times\cdots\times(a_n,b_n)$. Then
\[\partial^{(0,0,\tilde\alpha)}F(x)\approx_{\test'}\sum_{j_1=0}^{\alpha_1-1}\sum_{j_2=0}^{\alpha_2-1}
\partial^{\tilde\alpha}h_{1,2,j_1,j_2}(\tilde x_{1,2})x_1^{j_1}x_2^{j_2}.\]
Fixing now $c\in\ns\ster(a_1,b_1)$ and $d\in\ns\ster(a_2,b_2)$, we choose $\psi_m$ as before, $\phi(\tilde x_{1,2})\in\Fin(\ster{\tilde\Omega})$, multiply the previous expression by $\psi_m(c-x_1)\psi_m(d-x_2)\phi(\tilde x_{1,2})$ and
integrate over $\ster\Omega$ to obtain similarly that
\[\partial^{\tilde\alpha}F(c,d,\tilde x_{1,2}) \approx_{\test'} \sum_{j_1=0}^{\alpha_1-1} \sum_{j_2=0}^{\alpha_2-1} \partial^{\tilde\alpha} h_{1,2,j_1,j_2}(\tilde x_{1,2}) c^{j_1}d^{j_2}.\]
Now we substitute $c$ by $\alpha_1$ different values $c_1$, \dots, $c_{\alpha_1}$ $\in\ns\ster(a_1,b_1)$ and $d$ by $\alpha_2$ different values $d_1$, \dots, $d_{\alpha_2}$ $\in\ns\ster(a_2,b_2)$, with $c_i\not\approx c_j$ if $i\ne j$ and $d_i\not\approx d_j$ if $i\ne j$. The resulting linear system has $\alpha_1\alpha_2$ equations and $\alpha_1\alpha_2$ unknown functions $\partial^{\tilde\alpha} h_{1,2,j_1,j_2}$. The matrix of the system is (if the equations and unknowns are written down in a suitable order) the Kronecker-product (sometimes also called direct product, see e.g.~\cite{lancaster}) of the Vandermonde-matrices $(c_i^{j-1})_{i,j=1,\dots,\alpha_1}$ and $(d_i^{j-1})_{i,j=1,\dots,\alpha_2}$, with determinant \[\prod_{i<j}(c_j-c_i)^{\alpha_2}\prod_{i<j}(d_j-d_i)^{\alpha_1}\not\approx 0.\]
Another application of the previous lemma yields that
\[
h_{1,2,j_1,j_2}(\tilde x_{1,2}) \approx_{\test'}g_{1,2,j_1,j_2}(\tilde x_{1,2}) +\sum_{i=3}^n\sum_{k=0}^{\alpha_i-1}\tilde h_{ik}(\tilde x_{1,2,i})x_i^k\,,
\]
for some $g_{1,2,j_1,j_2}\in\ster{\Cnt\infty}(\Omega)$, S-continuous and finite-valued on $\ns(\ster\Omega)$ and $\tilde h_{ik}\in D'(\Omega)$. Substituting these expressions (for all $h_{i_1,i_2,j_1,j_2}$) in formula~\eqnref{induction-hyp} and absorbing the terms $g_{i_1 i_2 j_1 j_2}(\tilde x_{i_1 i_2})x_{i_1}^{j_1}x_{i_2}^{j_2}$ in the $g_{ij}(\tilde x_i)x_i^j$, we find formula~\eqnref{induction-result}.

Repeatedly applying this procedure, we conclude that
\[f(x)\approx_{\test'}\sum_{i=1}^n\sum_{j=0}^{\alpha_i-1} g_{ij}(\tilde x_i)x_i^j
+\sum_{j_1,j_2,\dots,j_n} c_{j_1,\dots,j_n}x_{1}^{j_1}x_{2}^{j_2}\cdots x_{n}^{j_n}\,,\]
for some $g_{ij}\in\ster{\Cnt\infty}(\Omega)$, S-continuous and finite-valued on $\ns(\ster\Omega)$ and constant $c_{j_1,\dots,j_n}\in D'(\Omega)$.
As a constant function belonging to $D'(\Omega)$ is necessarily
$\in\Fin(\ster\C)$ (see theorem \ref{pointvals}), we can absorb the terms
$c_{j_1,\dots,j_n}x_{1}^{j_1}x_{2}^{j_2}\cdots x_{n}^{j_n}$ in the $g_{ij}(\tilde x_i)x_i^j$ and finally obtain the required formula.
\end{proof}

Finally, before proving our second main result, we need a lemma of Robinson's \cite[Th.~5.3.14]{Robinson}. Robinson works with real-valued distributions on $\R$. We show that the result can be generalised to our situation.
\begin{lemma}\label{pointvals}
Let $T\in\test'(\Omega)$. If there exists a representative $f$ of $T$ which is S-continuous at $a\in\Omega$, then $f(a)\in\Fin(\ster\C)$.
Moreover, the value $\st f(a)$ does not depend on the chosen S-continuous
representative.
\end{lemma}
\begin{proof}
Let $\eps\in\R^+$. By S-continuity, there exists $r\in\R^+$ such that
$\abs{f(x)-f(a)}$ $\le \eps$, $\forall x\in\ster B(a,r)\subseteq\ster\Omega$.
Now let $\phi\in\test(B(a,r))$, real-valued,
$\phi(x)\ge 0$, $\forall x\in\Omega$ and $\int_\Omega\phi=1$. Then
\begin{align*}
&\abs{\int_{\ster\Omega} f(x)\ster\phi(x)\,dx-f(a)}\\
=&\abs{\int_{\ster B(a,r)}(f(x)-f(a))\ster\phi(x)\,dx}
\le\int_{\ster B(a,r)}\eps\abs{\ster\phi(x)}\,dx=\eps.
\end{align*}
As $f$ represents $T$, also $\abs{T(\phi)-f(a)}\le 2\eps$. In particular,
$f(a)\in\Fin(\ster\C)$. For any representative $g$ of $T$, S-continuous at $a$, we have the same inequality (possibly only for some smaller $r\in\R^+$), so $\abs{f(a)-g(a)}\le 4\eps$. As $\eps\in\R^+$ arbitrarily, $\st f(a)=\st g(a)$.
\end{proof}

\begin{theorem}
Let $f\in\ster{\Cnt\infty}(\Omega)$. Then $f\approx_{\test'(\Omega)}0$ iff
for each $K\csubset\Omega$, there exists $\alpha\in\N^n$ and
$g\in\ster\test(\Omega)$ such that $g(x)\approx 0$,
$\forall x\in \ster \Omega$ and $f=\partial^\alpha g$ on $\ster K$.
\end{theorem}
\begin{proof}
1. $\Rightarrow$: We first consider the case where $K\csubset\Omega$ is an interval.\\
Take an interval $K'\csubset\Omega$ with $K\csubset\inw{K'}$. By
theorem~\ref{local-distr}, there exists $h \in D'(\Omega)$ which is finite-valued and
S-continuous on $\ster K'$ and such that $\partial^\alpha h=f$ on $\ster K'$.
By lemma~\ref{Scont-prim} applied on the open interval $\tilde\Omega:=\inw{K'}$,
we find in particular that $h$ is $\test'(\tilde\Omega)$-infinitely close to
some $\tilde h\in\ster{\Cnt\infty}(\tilde\Omega)$, which is S-continuous on
$\ns(\ster {\tilde\Omega})$. As
$\tilde h(x)=\sum_{i=1}^n\sum_{j=0}^{\alpha_i-1} g_{ij}(\tilde x_i)x_i^j$, we
see that $\partial^\alpha\tilde h=0$ on $\ster{\tilde\Omega}$.
Now $h-\tilde h\approx_{\test'(\tilde\Omega)} 0$ and is S-continuous on $\ns(\ster
{\tilde\Omega})$, so by lemma~\ref{pointvals}, $h(x)-\tilde h(x)\approx 0$,
$\forall x\in \ns(\ster{\tilde\Omega})$. Further,
$\partial^\alpha (h-\tilde h)=\partial^\alpha h=f$ on $\ster K$. If
$\phi_0\in\test(\tilde\Omega)$ with $\phi_0=1$ on a neigbourhood of $\ster K$,
$g:=(h-\tilde h)\ster\phi_0$ has the required properties.\\
2. The general case, as well as the $\Leftarrow$-part follow in a way similar to theorem~\ref{local-distr}.
\end{proof}

\section{Hoskins and Sousa Pinto's theory of distributions}
In~\cite{hoskins}, R.~F.~Hoskins and J.~Sousa Pinto introduce another
nonstandard theory of distributions. In this setting, nonstandard
representatives of a distribution are {\it by definition} locally finite-order
derivatives of finite-valued and S-continuous functions. By
theorem~\ref{local-distr}, it now follows that representatives of distributions
in the sense of Hoskins and Sousa Pinto are exactly representatives of
distributions in the sense of Stroyan and Luxemburg.

\end{document}